\documentclass{amsart}
\usepackage{amsmath,amsfonts,amssymb,amsthm,amsxtra,latexsym,verbatim}

\newtheorem{theorem}{Theorem}[section]

\theoremstyle{definition}
\newtheorem{definition}[theorem]{Definition}

\theoremstyle{remark}

\numberwithin{equation}{section}

\newcommand{\N}{\mathbb{N}}
\newcommand{\R}{\mathbb{R}}
\newcommand{\Z}{\mathbb{Z}}

\newcommand{\M}{\mathcal{M}}

\renewcommand{\a}{\mathfrak{a}}
\newcommand{\p}{\mathfrak{p}}

\DeclareMathOperator{\Hom}{Hom}
\DeclareMathOperator{\Spec}{Spec}

\begin{document}

\title{Fan is to monoid as scheme is to ring}
\author{Howard~M Thompson}
\address{Department~of~Mathematics, University~of~Michigan, Ann~Arbor, Michigan 48109}
\email{hmthomps@umich.edu}
\date{}

\begin{abstract}
This paper generalizes the notion of a toric variety. In particular, these generalized toric varieties include examples of non-normal non-quasi-projective toric varieties. Such an example seems not to have been noted before in the literature.

This generalization is achieved by replacing a fan of strictly rational polyhedral cones in a lattice with a monoided space, that is a topological space equipped with a distinguished sheaf of monoids. For a classic toric variety, the underlying topological space of this monoided space is its orbit space under the action of the torus. And, if $\sigma$ is a cone, then the stalk of the structure sheaf of the monoided space at the point corresponding to $\sigma$ is the monoid $S_{\sigma}=\sigma\spcheck\cap M$ that is usually associated to $\sigma$. When applied to the affine toric variety associated to some cone $\sigma$, the monoided space so obtained is isomorphic to the spectrum of the monoid $S_{\sigma}$, where the spectrum of a monoid is defined in analogy with the definition of the spectrum of a ring. More generally, we will call a monoided space that is locally isomorphic to the spectrum of a monoid a fan and we form schemes from these fans by taking monoid algebras and glueing.
\end{abstract}

\maketitle

\section*{Introduction}

In recent years, various mathematicians have studied non-separated toric varieties or not necessarily normal toric varieties or combinatorial gadgets for describing toric varieties other than a fan of strictly rational polyhedral cones in a lattice. For examples, see A'Campo-Neuen \& Hausen~\cite{AH02}, Sturmfels~\cite{bS96} and Berchtold \& Hausen~\cite{BH04} respectively.

Following DeMeyer, Ford \& Miranda~\cite{DFM93}, we define a topology on a fan $\Delta$ in $\R^d$ by declaring the open sets to be its subfans. Then, like Kato~\cite{kK94}, we make such fans into monoided spaces by associating a sheaf of monoids to each $\Delta$. Our sheaf of monoids is similar to but differs from Kato's. Observing that this new monoided space is locally isomorphic to the spectrum of some monoid in the same sort of way a scheme is locally isomorphic to the spectrum of some ring, we define any monoided space that is locally isomorphic to the spectra of monoids to be a fan. The monoid algebra functor can then be used to associate a scheme to such a fan.

The advantages of the approach presented here are: 1) some of the intuition built up from studying schemes can be carried over to this new gadget; 2) for any toric variety, not just the projective ones, the gadget presented here comes equipped with an abstract analog of the moment map and often torus invariant objects on the toric variety descend to this gadget via this map; 3) when the toric variety is projective this new gadget is the corresponding polytope equipped with some extra structure and the associated abstract moment map is, essentially, the moment map; 4) this gadget permits the consideration of non-normal toric varieties that are neither affine nor projective, in fact the scheme under consideration need not be a variety; 5) at the same time, the collection of schemes built using these gadgets is small enough that all the normal varieties in it are classic toric varieties.

\section{Basic Notions}

\begin{definition}
Let $S$ be a commutative monoid. We say a subset $\a\subseteq S$ is an \emph{ideal} if $\a+S\subseteq\a$. We say an ideal $\p\subset S$ is \emph{prime} if its complement is a submonoid of $S$.

The \emph{spectrum} of $S$ is a pair $(\Spec S,\M)$ consisting of a topological space $\Spec S$ and a sheaf of monoids $\M$ on $\Spec S$ defined as follows: The underlying set of $\Spec S$ is the collection of prime ideals of $S$, the open sets of the form $D(f)=\{\p\in\Spec S\mid f\in S\setminus\p\}$ form a base for the topology of $\Spec S$, and $\Gamma(D(f),\M)=S+\N(-f)$. I should warn you that I will use the same notion for both rings and monoids. In this paper, if $A$ is a ring, then $\Spec A$ is the set of prime ideals of $A$ as a ring with its standard topology, etc.

A \emph{monoided space} is a pair $(\Delta, \M)$ consisting of a topological space $\Delta$ and a sheaf of monoids $\M$ on $\Delta$. A morphism of monoided spaces consists of a continuous map $\varphi:\Delta\to\Delta'$ and local homomorphism of sheaves of monoids $\varphi^{\#}:\M_{\Delta'}\to\varphi_{*}\M_{\Delta}$. That is, for every point $\sigma\in\Delta$, the homomorphism of monoids $\varphi_{\sigma}^{\#}:\M_{\Delta',\varphi(\sigma)}\to\varphi_{*}\M_{\Delta,\sigma}$ maps non-units to non-units.

We say a monoided space $(\Delta,\M)$ is a \emph{fan} if every point of $\Delta$ has an open neighborhood $U$ such that there exists a monoid $S$ with $(U,\M|_U)\cong(\Spec S,\M_{\Spec S})$.
\end{definition}

Every monoid $S$ has a unique maximal ideal $S^+=S\setminus S^*$ so there is no distinct notion of a locally monoided space. Furthermore, every point $\sigma$ on a fan $(\Delta,\M)$ has a unique smallest open neighborhood and this neighborhood is isomorphic to $\Spec\M_{\sigma}$.

If $\Delta$ is a fan and $A$ is a ring, we associate a scheme to this pair using the monoid algebra functor $A\mapsto A[\Delta]$: That is, if $\Spec S$ and $\Spec S'$ are open subsets of the fan $\Delta$, and $f\in S$ is such that $D(f)\subseteq\Spec S'$, we glue the affine schemes $\Spec A[S]$ and $\Spec A[S']$ along $\Spec S_f$. Here the map $A[S_f]\to A[S']$ is the one induced by the restriction map $S_f\to S'$ from $\Delta$. This makes sense because the canonical inclusions of the form $S\hookrightarrow A[S];\,s\mapsto t^s$ (We write elements of $A[S]$ or more expansively, $A[t;S]$, as polynomials in $t$ with coefficients in $A$ and exponents in $S$.) induce a continuous map $A[\Delta]\to\Delta$. The pullback of $D(s)\subseteq\Spec S$ is $D(t^s)\subseteq\Spec A[S]$. More generally, if $\Delta$ is a fan and $X$ is a scheme, we may associate the $X$-scheme $X[\Delta]=X\times_{\Spec\Z}\Z[\Delta]\xrightarrow{pr_1}X$ to this pair.

\section{The fan of a classic toric variety}

\begin{definition}
Let $\Delta$ be a fan in a lattice $N\cong\Z^d$, that is, let $\Delta$ be a collection of strongly convex rational polyhedral cones in the $\R$-vector space $N_{\R}=N\otimes_{\Z}\R$. We associate to $\Delta$ a monoided space as follows. The underlying topological space, $\Delta_{top}$, is the collection $\{\sigma\in\Delta\}$ equipped with the topology whose open sets are the subfans of $\Delta$. That is, the open sets are subsets of $\Delta$ that are also fans in $N$. And, the sheaf of monoids on $\Delta_{top}$ is the sheaf $\M$ such that $\Gamma(\sigma_{top},\M|_{\sigma})=S_{\sigma}(=\sigma\spcheck\cap M)$. Here we have identified the cone $\sigma$ with the fan in $N$ consisting of $\sigma$ and its faces and $M=\Hom_{\Z}(N,\Z)$.
\end{definition}

\begin{theorem}
If $\Delta$ is a fan in $N$, then $(\Delta_{top},\M)$ is a fan.
\end{theorem}

\begin{proof}
It suffices to prove $(\sigma_{top},\M|_{\sigma})\cong\Spec S_{\sigma}$ for all $\sigma\in\Delta$. This is straightforward.
\end{proof}

Notice that, if one starts with the normal fan of a polytope, then topological space we obtain is isomorphic to the set of faces of the polytope equipped with the topology such the the star of a face is the smallest open subset containing that face. More generally, it is worth noting that $\Delta_{top}$ is homeomorphic to the orbit space of the classic toric variety associated to $\Delta$.

\section{When is $k[\Delta]$ a normal $k$-variety?}

This section is devoted to proving: If $\Delta$ is a fan and $k$ is a field,  then $k[\Delta]$ is a normal $k$-variety if and only if it is a classic toric variety. Here by $k$-variety we mean an integral, separated scheme of finite type over $k$.

\begin{theorem}\cite[Exercise~II.3.5]{rH77}
If $f:X\to Y$ is a morphism of schemes of finite type, then for every open affine subset $V=\Spec B\subseteq Y$, and for every open affine subset $U=\Spec A\subseteq f^{-1}(V)$, $A$ is a finitely generated $B$-algebra.
\end{theorem}

In particular, if $k[\Delta]$ is a $k$-scheme of finite type, then the stalks of $\M$ are all finitely generated monoids. Recall that a monoid $S$ is said to be torsion-free if $ns=ns'$ for some positive integer $n$ implies $s=s'$ in $S$.

\begin{theorem}\cite[Theorem~8.1]{rG84}
Let $A$ be a nonzero ring. The monoid algebra $A[S]$ is an integral domain if and only if $A$ is an integral domain and $S$ is torsion-free and cancellative.
\end{theorem}

As a consequence of this proposition, if $k[\Delta]$ is integral, then the stalks of $\M$ are all cancellative, torsion-free monoids. In particular, if $k[S]$ is a domain, the $S$ embeds in a group it generates. Use $\Z S$ to denote this group and identify $S$ with its image in $\Z S$.

If $S$ is a finitely generated, cancellative torsion-free monoid $\Z S$ is a finitely generated, torsion-free Abelian group. That is, $\Z S\cong\Z^n$ for some non-negative integer $n$. Consider the rational polyhedral cone $\R_{\geq0}S=\left\{\sum_{i=1}^l r_i\otimes s_i\mid r_i\geq0,\,s_i\in S\right\}$ in the finitely generated $\R$-vector space, $\R S=\R\otimes_{\Z}\Z S$, generated by $S$. Identify $S$ and $\Z S$ with their images in $\R S$. If $s\in\Z S\cap\R_{\geq0}S$, then $ns=s'\in S$ for some positive integer $n$ and some $s'\in S$. So, the image of $s$ in $k[\Z S\cap\R_{\geq0}S]$, $t^s$, satisfies the monic polynomial $x^n-t^{s'}$ with coefficients in $k[S]$. Therefore, if $s\in\Z S\cap\R_{\geq0}S$ and $k[S]$ is integrally closed, then $s\in S$. In other words, if $\Spec k[S]$ is a normal $k$-variety, then $S$ is a rational polyhedral cone in some finite dimensional vector space intersected with a sublattice generated by a basis of the vector space. More briefly, if $X=\Spec k[S]$ is a normal $k$-variety, then $X$ is an affine toric variety.

So far, we know that if $X=k[\Delta]$ is a normal integral scheme of finite type over $k$, then $X$ can be glued together from affine toric varieties along open torus invariant subvarieties. We now wish to show that if in addition $X$ is separated, then $X$ is a classic toric variety.

\begin{theorem}
Suppose $X=k[\Delta]$ is a normal integral scheme of finite type over $k$. If $X$ is separated, then $X$ is a classic toric variety.
\end{theorem}

\begin{proof}
We have already established that $X$ is glued together from a collection of affine toric varieties. We need to show that the cones corresponding to these toric varieties fit together to form a fan of strictly rational polyhedral cones in a lattice.

Since on each of the fans of these affine toric varieties the generic point is the unique smallest open set, all these generic points are identified by the gluing. So, there is a fixed lattice $M$, the stalk of $\M$ at the generic point, containing every monoid $S$ that occurs. Furthermore, $M=\Z S$ for every monoid $S$ and $\R_{\geq0}S$ is dual to a strictly convex rational polyhedral cone in $N$, where $N=\Hom_{\Z}(M,\Z)$. So, it makes sense to ask whether $X$ comes from a fan.

All we have to prove is: The intersection of any two cones in this collection is a face of each. We know the intersection $U=\Spec k[S_1]\cap\Spec k[S_2]\subseteq X$ is affine since $X$ is separated. Furthermore, $U$ is torus invariant. So, $U$ comes from a subfan of the cone, $\sigma$, in $N$ corresponding to $\Spec k[S_1]$. Since $U$ is also affine, $U$ corresponds to a face of $\sigma$.
\end{proof}

\bibliographystyle{amsplain}

\providecommand{\bysame}{\leavevmode\hbox to3em{\hrulefill}\thinspace}
\providecommand{\MR}{\relax\ifhmode\unskip\space\fi MR }
\providecommand{\MRhref}[2]{%
  \href{http://www.ams.org/mathscinet-getitem?mr=#1}{#2}
}
\providecommand{\href}[2]{#2}


\end{document}